\documentclass[reqno,10pt,centertags]{amsart}
\usepackage{amsmath,amsthm,amscd,amssymb,latexsym,upref}

\newcommand{\bbN}{{\mathbb{N}}}
\newcommand{\bbR}{{\mathbb{R}}}

\newcommand{\bbZ}{{\mathbb{Z}}}
\newcommand{\bbC}{{\mathbb{C}}}

\newcommand{\calB}{{\mathcal B}}

\newcommand{\calH}{{\mathcal H}}

\newcommand{\calU}{{\mathcal U}}

\newcommand{\no}{\nonumber}
\newcommand{\lb}{\label}
\newcommand{\f}{\frac}
\newcommand{\ul}{\underline}
\newcommand{\ol}{\overline}
\newcommand{\ti}{\tilde  }
\newcommand{\wti}{\widetilde  }

\newcommand{\spec}{\text{\rm{spec}}}

\newcommand{\dom}{\text{\rm{dom}}}

\newcommand{\supp}{\text{\rm{supp}}}

\newcommand{\bi}{\bibitem}

\newcommand{\beq}{\begin{equation}}
\newcommand{\eeq}{\end{equation}}
\newcommand{\ba}{\begin{align}}
\newcommand{\ea}{\end{align}}

\renewcommand{\Re}{\text{\rm Re}}
\renewcommand{\Im}{\text{\rm Im}}
\renewcommand{\ln}{\text{\rm ln}}

\allowdisplaybreaks
\numberwithin{equation}{section}

\newtheorem{theorem}{Theorem}[section]
\newtheorem{lemma}[theorem]{Lemma}

\theoremstyle{definition}

\theoremstyle{remark}
\newtheorem{remark}[theorem]{Remark}

\begin{document}
\title[An inverse problem for point inhomogeneities]
{An Inverse Problem for Point Inhomogeneities}
\author[Gesztesy and Ramm]{Fritz Gesztesy and 
Alexander G.~Ramm}
\address{Department of Mathematics,
University of
Missouri, Columbia, MO
65211, USA}
\email{fritz@math.missouri.edu\newline
\indent{\it URL:}
http://www.math.missouri.edu/people/fgesztesy.html}
\address{Department of Mathematics, Kansas State University,
Manhattan, KS 66506-2602, USA}
\email{ramm@math.ksu.edu}
\thanks{This paper was written during AGR's visit at the 
Mathematics Department of the University of Missouri, 
Columbia, in connection with his Big 12 Fellowship.}

\begin{abstract}
We study quantum scattering theory off $n$ 
point inhomogeneities ($n\in\bbN$) in three dimensions. The
inhomogeneities  (or generalized point interactions) 
positioned at $\{\xi_1,\dots,\xi_n\}\subset\bbR^3$ are 
modeled  in terms  of the $n^2$ (real) parameter family of
self-adjoint extensions of 
$-\Delta\big|_{C^\infty_0(\bbR^3\backslash
\{\xi_1,\dots,\xi_n\})}$ in $L^2(\bbR^3)$. The Green's 
function, the scattering 
solutions and the scattering amplitude for this model are
explicitly  computed in terms of elementary functions. 
Moreover, using the  connection between fixed energy 
quantum scattering and  acoustical  scattering, the 
following inverse spectral result in  acoustics is proved: 
The knowledge of the scattered field  on a plane outside 
these point-like inhomogeneities, with  all inhomogeneities
located  on one side of the plane,  uniquely determines the 
positions and boundary conditions associated with them. 
\end{abstract}

\maketitle

\section{Introduction}\lb{s1}

To describe the inverse problem solved in this paper in 
some detail we need a few preparations. Let $\bbR_\pm^3 
=\{x=(x_1,x_2,x_3)\in\bbR^3 \,|\, x_3 \gtrless 0\}$, 
$P=\{x=(x_1,x_2,x_3)\in\bbR^3 \,|\, x_3 =0\}$, 
$D\subset\bbR^3_-$ a domain with smooth boundary  
and finitely many connected components, 
and $v\in L^2(D)$, $\supp (v)\subseteq D$, $v$
real-valued . Consider the fixed energy scattering
problem
\begin{equation}
\begin{cases}
[-\bigtriangledown^2_x-k_0^2-k_0^2
v(x)]u(x,y)=\delta(x-y), \quad x,y\in\bbR^3, \, x\neq y,\\
\lim_{|x|\to\infty}|x|\big[\f{\partial}{\partial |x|}
u(x,y) -ik_0 u(x,y)\big]=0 \quad \text{uniformly in  
directions $\omega=|x|^{-1}x$}\\
\hspace*{4.4cm} \text{and uniformly in $y$ for $y$ varying
in compact sets,}
\end{cases} \lb{1.1}
\end{equation}
with $k_0>0$ (the wave number) a fixed positive constant. 
Here $c(x)=[v(x)+1]^{-1/2}$ has the
physical meaning of the  wave velocity profile in the 
medium, $v(x)$ is
the inhomogeneity in the velocity profile, and $u(x,y)$
represents the acoustic pressure generated by a point
source at the point $y\in\bbR^3$.

The inverse problem (IP) associated with \eqref{1.1}, more
precisely, the inversion of the surface data $u(x,y)$ for the
velocity profile $c(x)$, then can be formulated as follows:\\

\noindent {\bf IP~1.1.}~Given the data 
$\{u(x,y)\}_{\substack{x,y\in P\\x\neq y}}$ at fixed
$k_0>0$, determine 
$v(x)$, $x\in D$.\\

A solution of this inverse problem (i.e., uniqueness of
$v(x)$ and recovery of $v(x)$ {}from the prescribed data) is
described in \cite[Sects.~III.6, IV.2]{Ra92}. Numerical
methods in connection with IP~1.1 are discussed in
\cite[Sect.~V.3]{Ra92}. 

Since $u(x,y)$ in \eqref{1.1} can be identified with the
Green's function at fixed energy $k_0^2>0$, 
\begin{equation}
G(k_0^2,x,y)=(-\Delta -k_0^2v-k_0^2)^{-1}(x,y), \quad 
x,y\in\bbR^3, \, x\neq y, \lb{1.2}
\end{equation}
associated with the self-adjoint (Schr\"odinger-type) operator 
\begin{equation}
H=-\Delta -k_0^2v, \quad \dom(H)=H^{2,2}(\bbR^3) \lb{1.3} 
\end{equation}
in $L^2 (\bbR^3)$, we can reformulate the inverse problem
IP~1.1 in the following equivalent form: \\

\noindent {\bf IP~1.1'.}~Given the data 
$\{G(k_0^2,x,y)\}_{\substack{x,y\in P\\x\neq y}}$ at a fixed
energy
$k_0^2>0$, determine $v(x)$, $x\in D$. \\

For practical applications in connection with ultrasound
mammography tests (as opposed to x-ray mammography) and in
the area of material science in connection with the
detection of cracks and cavities, it is of relevance to
consider inhomogeneities $v(x)$ of the special form 
\begin{equation}
v(x)=\sum_{j=1}^n v_j(x-\xi_j), \quad 
\{\xi_1,\dots,\xi_n\}\subset\bbR^3_-, \lb{1.4}
\end{equation}
where $v_j\in L^2 (D_j)$, $\supp(v_j(\cdot-\xi_j))\subseteq 
D_j$, and $D_j\subset\bbR^3_-$ are connected domains
with smooth boundaries and sufficiently small 
diameters $d_j$ with respect
to the wave length (i.e., $(\max_{1\leq j\leq n}d_j)k_0
\ll 1$).~A numerical procedure recovering the $\xi_j$ (and 
hence the approximate position of the small inhomogeneities)  
and the intensities of the inhomogeneities, defined by $V_j=
\int_{D_j} dx\, v_j(x-\xi_j)$, $1\leq j\leq n$, has recently
been discussed in \cite{Ra99}.

At this point we are in a position to describe the inverse
problem considered in this paper. In view of the physical
applications mentioned in connection with \eqref{1.4}, we
will now consider the idealized situation
of inhomogeneities $v_j(x-\xi_j)$ of {\it point-like}
support at $\xi_j$, $1\leq j\leq n$. Intuitively, we 
want to solve the inverse problem 
\begin{equation}
\begin{cases}
[-\bigtriangledown^2_x-k_0^2-k_0^2\sum_{j=1}^n 
v_j(x-\xi_j)]u(x,y)=\delta(x-y), \quad x,y\in\bbR^3, \, 
x\neq y,\\
\lim_{|x|\to\infty}|x|\big[\f{\partial}{\partial |x|}
u(x,y) -ik_0 u(x,y)\big]=0 \quad \text{ uniformly in  
directions $\omega=|x|^{-1}x$}\\
\hspace*{4.5cm} \text{and uniformly in $y$ for $y$ varying
in compact sets,}
\end{cases} \lb{1.5}
\end{equation}
where formally
\begin{equation}
v_j(x-\xi_j)=a_j\delta(x-\xi_j), \quad 1\leq j\leq n, 
\lb{1.6}
\end{equation}
for some ``coupling'' constants $a_j\in\bbR$, $1\leq j \leq
n$. However, as is well-known, point-like inhomogeneities
of the type \eqref{1.6} as potential 
coefficients in a Schr\"odinger-type operator in dimensions 
$d\geq 2$ do not lead to an operator or quadratic form 
perturbation of the Laplacian $-\Delta$, $\dom(-\Delta)
=H^{2,2}(\bbR^d)$ in $L^2(\bbR^d)$, where $d$ denotes the 
corresponding space dimension. One possible
way around  this difficulty for $d=2$ and $d=3$ is the introduction of an
appropriate coupling constant renormalization procedure. This 
point of
view is presented in detail in \cite[Ch.~II.1]{AGHKH88}. 
Alternatively to
this renormalization procedure for $d=2,3$, one can apply the 
theory of
self-adjoint extensions of closed symmetric densely defined
linear operators in a Hilbert space to the operator
\begin{equation}
\ol{-\Delta\big|_{C_0^\infty (\bbR^3\backslash\Xi)}}, 
\quad \Xi=\{\xi_1,\dots,\xi_n\}\subset\bbR^3_- \lb{1.7}
\end{equation}
in $L^2 (\bbR^3)$. (Here the $\ol T$ denotes the operator 
closure of $T$ and we refer to Remark~\ref{r3.5} for a brief 
discussion of the situation in different dimensions 
$d\in\bbN$.)~In this
paper we follow the latter approach  and model the Laplacian
$-\Delta$ perturbed by point-like perturbations of the type
$-k_0^2\sum_{j=1}^n a_j\delta(x-\xi_j)$ by self-adjoint 
extensions of 
$\ol{-\Delta\big|_{C_0^\infty (\bbR^3\backslash\Xi)}}$,
denoted by $-\Delta_{\theta,\Xi}$, parametrized by the 
$n^2$ (real) parameter family of self-adjoint matrices
$\theta$ in $\bbC^n$. 

Taking advantage of the equivalence of the inverse problems
IP~1.1 and IP~1.1', we can now formulate the
inverse problem associated with point-like inhomogeneities,
as studied in this paper, in a precise manner as follows: \\

\noindent {\bf IP~1.2.}~Prove that the data
$\{G_{\theta,\Xi}(k_0^2,x,y)\}_{\substack{x,y\in P\\x\neq
y}}$  at fixed energy $k_0^2>0$, uniquely determine
$\Xi=\{\xi_1,\dots,\xi_n\}\subset\bbR^3_-$ and the
self-adjoint $n\times n$ matrix $\theta$ in $\bbC^n$.
\\

Here $G_{\theta,\Xi}(k_0^2,x,y)$ denotes the Green's 
function associated with $-\Delta_{\theta,\Xi}$, that is, 
\begin{equation}
G_{\theta,\Xi}(z,x,y)=(-\Delta_{\theta,\Xi}-z)^{-1}(x,y), 
\quad \det(P_{\theta,\Xi}(z))\neq 0, \,\,
x,y\in\bbR^3\backslash\Xi, 
\, x\neq y. \lb{1.8}
\end{equation}

While IP~1.1 (resp.,
(IP~1.1') is concerned with uniqueness and reconstruction
of $v(x)$,
$x\in\ol D$, IP~1.2, as studied in this paper, focuses on the
unique determination of $\Xi$ and $\theta$ by the data
measured on the plane $P$. 

In Section~\ref{s2} we present a detailed account of
Krein's formula of self-adjoint extensions of closed
symmetric operators in a Hilbert space, our principal tool
in describing the $n^2$ (real) parameter family of
self-adjoint extensions $-\Delta_{\theta,\Xi}$ of
\eqref{1.7} in Section~\ref{s3}. In particular, we
explicitly describe the Green's function, the scattering
solutions, and the scattering amplitude associated with 
$-\Delta_{\theta,\Xi}$ in Section~\ref{s3}. The inverse
problem IP~1.2 is solved in our final Section~\ref{s4}. 

\section{Krein's Formula for Self-adjoint Extensions} \lb{s2}

In this section we recall Krein's formula,
which describes the resolvent difference of two self-adjoint
extensions $A_1$ and $A_2$
of a densely defined closed symmetric linear operator
$A$ with
deficiency indices $(n,n)$, $n\in \bbN$. (Reference
\cite{GMT98} treats this topic in the general case where
$n\in\bbN\cup\{\infty\}$. Here we
restrict ourselves to the case $n<\infty$.) We start with 
the basic setup following \cite{AG93}.

Let $\calH$ be a separable complex Hilbert space,
$\dot A:\dom(\dot A)\to\calH, \,\, \ol{\dom(\dot A)}
=\calH$ a
densely defined closed symmetric linear operator in
$\calH$ with finite
and equal deficiency indices
 ${\rm def}(\dot A)=(r,r)$,
$r\in \bbN$. Let $A_\ell, \, \ell=1,2$,  be two distinct
self-adjoint
extensions
of $\dot A$  and denote by $A$ the maximal common part
of $A_1$ and $A_2$,
that is, $A$ is the largest closed extension of
$\dot A$  with
$\dom(A)=\dom(A_1)\cap \dom(A_2)$. In this case one calls 
$A_1$ and $A_2$ {\it relatively prime} with respect to $A$. 
Let $0\le p\le r-1$ be the
maximal number of elements in $\dom(A)=\dom(A_1)
\cap \dom(A_2)$
which are linearly independent modulo $\dom(\dot A)$.
Then $A$ has
 deficiency indices ${\rm def} (A)=(n,n)$, $n=r-p$.
Next, denote by
 $\ker (A^*-z), \, z\in \bbC\backslash\bbR$ the deficiency
subspaces
of $A$ and define
\begin{equation}\label{e1}
W_{1,z,z_0}=I+(z-z_0)(A_1-z)^{-1}=(A_1-z_0)(A_1-z)^{-1},
\,\,z,z_0\in
\rho (A_1),
\end{equation}
where $I$ denotes the identity operator in $\calH$ and
$\rho(T)$
abbreviates the resolvent set of $T$. One verifies
\begin{equation}\label{e2}
W_{1,z_0,z_1}W_{1,z_1,z_2}=W_{1,z_0,z_2},\,\, z_0, z_1,
z_2\in \rho (A_1)
\end{equation}
and
\begin{equation}\label{e3}
W_{1,z,z_0}\ker(A^*-z_0)=\ker(A^*-z).
\end{equation}
Let $\{ u_j(i) \}_{1\le j \le n}$ be an orthonormal basis
for
$ \ker(A^*-i)$ and define
\begin{equation}\label{e4}
u_{j}(z)=W_{1,z,i}u_{j}(i)=(A_1-i)(A_1-z)^{-1}u_j(i), 
\quad 1\le j\le n, \,\, z\in \rho(A_1).
\end{equation}
Then $\{ u_{j}(z) \}_{1\le j \le n}$ is a basis for
$\ker(A^*-z)$,
 $z\in \rho(A_1)$ and since
$ W_{1,-i,i}=(A_1-i)(A_1+i)^{-1}$  is the unitary
Cayley transform of $A_1$,
$\{ u_{j}(-i) \}_{1\le j \le n}$ is in fact an
orthonormal basis for $\ker(A^*+i)$.

The basic result on Krein's formula, as presented by
Akhiezer and
 Glazman \cite{AG93}, Sect. 84, then reads as follows.

\begin{theorem} \mbox{\rm (Krein's formula, 
\cite[Sect.~84]{AG93}.)} \lb{t1} \\
There exists a $P_{1,2}(z)=\big(P_{1,2}(z)_{j,j'}
\big)_{1\le j,j' \le n}\in M_n(\bbC)$,
$z\in \rho(A_2)\cap\rho(A_1)$, such that
\begin{align}
&\det (P_{1,2}(z))\ne 0, \,\,z\in \rho(A_2)\cap\rho(A_1),
\label{e5} \\
&P_{1,2}(z)^{-1}=P_{1,2}(z_0)^{-1}-(z-z_0)\big((u_{j}
(\bar z),u_{j'}( z_0))\big)_{1\leq j,j'\leq n}\, ,
\,\,z,z_0\in \rho(A_1), \label{e6} \\
&\Im \, (P_{1,2}(i)^{-1})=-I_n, \label{e7} \\
&(A_2-z)^{-1}=(A_1-z)^{-1}+
\sum_{j,j'=1}^nP_{1,2}(z)_{j,j'}(u_{j'} (\bar z), 
\cdot \,) u_{j} (z), \quad z\in \rho(A_2)\cap\rho(A_1). 
\label{e8}
\end{align}
\end{theorem}

Here $\Im(T)=(T-T^*)/(2i)$ and $\Re(T)=(T+T^*)/2$ denote 
the imaginary and real parts of the matrix $T$, 
respectively.

We note that $P_{1,2}(z)^{-1}$ extends by continuity
{}from $z\in\rho(A_2)\cap\rho(A_1)$ to all of $\rho(A_1)$
since the right-hand
side of
(\ref{e6}) is continuous for $z\in \rho(A_1)$. The
normalization condition
(\ref{e7}) is not mentioned  in \cite{AG93} but
it trivially follows
{}from (\ref{e6}) and the fact
\begin{equation}
(u_j(i), u_{j'}(i))=\delta_{j,j'}, \quad 1\le j,j' 
\le n \lb{e8a}
\end{equation}
(where $\delta_{j,j'}$ denotes Kronecker's symbol)
and {}from
\begin{equation}
P_{1,2}^*(z)=P_{1,2}(\bar z), \quad z\in\rho(A_1)
\cap\rho(A_2). \lb{e8b}
\end{equation}
Taking $z=\overline{z_0}$ in (\ref{e6}) shows that
$-P_{1,2}(z)^{-1}$
and hence
 $P_{1,2}(z)$ is a matrix-valued Herglotz function,
that is,
\begin{equation}\label{e8c}
\Im\, (P_{1,2}(z))>0, \,\, z\in \bbC_+.
\end{equation}
Strict positive
definiteness in (\ref{e8c}) follows {}from the fact that
$\{u_{j}(z)\}_{1\le j \le n}$ are linearly independent
for $z\in \bbC_+$
and hence
$((u_{j}(z), u_{j'}(z))_{1\le j,j' \le n} >0$.

Next we turn to the connection between $P_{1,2}(z)$ and
von Neumann's parametrization of self-adjoint extensions 
of $A$ as discussed in detail in \cite{GMT98}. Due
to (\ref{e6}), $P_{1,2}(z)^{-1}$ is determined for all
$z\in \rho(A_1)$
in terms of
$P_{1,2}(i)^{-1}, \,\,(A_1-z)^{-1}$ and
$\{ u_j(i)\}_{1 \le j \le n} $,
\begin{align}
&P_{1,2}(z)^{-1}=P_{1,2}(i)^{-1}-(z-i)I_n-(1+z^2)\big((u_j(i),
(A_1-z)^{-1}u_{j'}(i))\big)_{1 \le j,j' \le n}, \no \\
& \hspace*{9.5cm} z\in \rho(A_1). \label{e9}
\end{align}
Hence it suffices to focus on
\begin{equation}
P_{1,2}(i)^{-1}=\Re\,(P_{1,2}(i)^{-1})-iI_n. \lb{e9a}
\end{equation}
Let
\begin{equation}\label{e10}
\calU_\ell:\ker (A^*-i) \to \ker (A^*+i),
\quad \ell=1,2,
\end{equation}
be the linear isometric isomorphisms that parameterize
$A_\ell$ according
to von Neumann's formula
\begin{align}
&A_\ell(f+(I+\calU_\ell)u_+)=Af+i(I-\calU_\ell)u_+, \no \\
&\dom(A_\ell)=\{ (g+(I+\calU_\ell)u_+)\in \dom(A^*) \,
\vert \, g\in \dom(A), \, u_+\in \ker (A^*-i) \}, \no \\
&\hspace*{9.8cm} \ell=1,2. \lb{e11} 
\end{align}
Next, denote by
$U_\ell=(U_{\ell,j,j'})_{1\le j,j'
\le n}\in M_n(\bbC),\,\,  \ell=1,2$
the unitary matrix representation of $\calU_\ell$ with
respect to the bases
$\{ u_j(i) \}_{1 \le j \le n}$ and
$\{ u_{1,j}(-i) \}_{1 \le j \le n}$
of $\ker (A^*-i)$ and $\ker (A^*+i)$ respectively, that is,
\begin{equation}\label{e12}
\calU_\ell u_j(i)=\sum_{j'=1}^n U_{\ell,j',j}u_{1,j'}(-i),
\,\, 1\le j \le n, \,\,
\ell=1,2.
\end{equation}

\begin{lemma} \rm{(}\cite{GMT98}.\rm{)} \lb{l2}\\
\noindent \rm{(}i) $U_1=-I_n$. \\
\noindent \rm{(}ii\rm{)} $-1\notin\spec(U_2)$. \\
\noindent \rm{(}iii\rm{)} $U_\ell, \,\,\ell=1,2$ and $P_{1,2}(i)$
are  connected by
\begin{equation}\label{e19}
P_{1,2}(i)=\frac{i}{2} (I_n+U_2^{-1})=\frac{i}{2}(U_2^{-1}
-U_1^{-1}).
\end{equation}
\end{lemma}

Here $\spec(T)$ denotes the spectrum of $T$. 

Next, writing
\begin{equation}
U_2=\exp(i\theta_2), \quad \theta_2^*=\theta_2 \lb{2.20}
\end{equation}
for the matrix representation of $\calU_2$ with respect to the
bases $\{ u_j(i) \}_{1 \le j \le n}$ and 
$\{ u_{1,j}(-i) \}_{1 \le j \le n}$ of 
$\ker (A^*-i)$ and $\ker(A^*+i)$, one verifies 
\begin{equation}
\Re(P_{1,2}(i)^{-1})=\tan(\theta_2/2). \lb{2.22}
\end{equation}
Introducing the matrix-valued Herglotz function 
$M_1(z)$ associated with $A_1$ (cf. \cite{Do65}, \cite{GMT98}) 
by
\begin{equation}\label{e21}
M_1(z)=zI_n +(1+z^2) \big((u_j(i), (A_1-z)^{-1}
u_{j'}(i))_{1\le j,j' \le n}\big), \quad z\in \rho(A_1),
\end{equation}
$P_{1,2}(z)$ in Krein's formula \eqref{e8} then can be 
rewritten as
\begin{align}
P_{1,2}(z)&=(\tan(\theta_2/2)-M_1(z))^{-1} \no \\
&=\big(\tan(\theta_2/2)-zI_n-(1+z^2)\big((u_j(i), (A_1-z)^{-1}
u_{j'}(i))_{1\le j,j' \le n}\big)\big)^{-1},  \no \\
&\hspace*{8cm} \quad z\in \rho(A_1). \lb{2.23}
\end{align}
We emphasize that 
\begin{equation}
\{((1/2) +m)\pi\}_{m\in\bbZ}\notin\spec(\theta_2) \lb{2.23b}
\end{equation}
according to Lemma~\ref{l2}\,(ii), due to our hypothesis that 
$A_1$ and $A_2$ are relatively prime with respect to $A$.

For subsequent purposes it is useful to introduce the 
self-adjoint operator $\vartheta_2\in\calB(\ker(A^*-i))$ 
defined 
through its matrix representation $\theta_2$ with respect 
to the basis $\{ u_j(i) \}_{1 \le j \le n}$, that is,
\begin{equation}
\theta_{2,j,j'}=(u_j(i),\vartheta_2 u_{j'}(i)), \quad 
1\leq j,j'\leq n. \lb{2.23a}
\end{equation}

The discussion of Krein's formula thus far dealt 
exclusively with the 
orthonormal bases $\{ u_j(i) \}_{1 \le j \le n}$ and 
$\{ u_{1,j}(-i) \}_{1 \le j \le n}$ of 
$\ker (A^*-i)$ and $\ker(A^*+i)$ following our discussion 
in \cite{GMT98} and \cite[Appendix~B]{GT97}. In the
remainder of this paper, however, it will be be more convenient 
to discuss matrix representations of $M_1(z)$ and $\calU_2$ 
with respect to a natural (cf.~the comment following \eqref{3.6}), 
but not
necessarily orthogonal  basis. Hence we briefly discuss the 
effect of a
change of  basis in connection with Krein's formula \eqref{e8}. 
Let 
$\{\ti u_j(i)\}_{1\leq j \leq n}$ be another (not necessarily 
orthogonal basis) of $\ker(A^*-i)$ and define
\begin{align}
&\ti u_{j}(z)=(A_1-i)(A_1-z)^{-1}\ti u_j(i), 
\quad 1\le j\le n, \,\, z\in \rho(A_1), \label{2.24} \\
&\calU_\ell \ti u_j(i)=\sum_{j'=1}^n 
\wti U_{\ell,j',j}\ti u_{j'}(-i), \quad 1\leq j\leq n, \,\, 
\ell=1,2, \lb{2.25} \\
&\wti U_2=\exp(i\ti \theta_2), \quad 
{\ti \theta_2}^*=\ti\theta_2. \lb{2.26}
\end{align}
In addition, one verifies
\begin{equation}
\wti U_1=-I_n \lb{2.27}
\end{equation}
as in Lemma~\ref{l2}\,(i). Krein's formula \eqref{e8} then can 
be rewritten in the form
\begin{equation}
(A_2-z)^{-1}=(A_1-z)^{-1}+
\sum_{j,j'=1}^n \wti P_{1,2}(z)_{j,j'}(\ti u_{j'} (\bar z), 
\cdot \,) \ti u_{j} (z), \quad z\in \rho(A_2)\cap\rho(A_1), 
\label{2.28}
\end{equation}
where
\begin{align}
\wti P_{1,2}(z)&=(\tan(\ti \theta_2/2)-\wti M_1(z))^{-1} \no \\
&=\big(\tan(\ti \theta_2/2)-zI_n-(1+z^2)\big((\ti u_j(i),
(A_1-z)^{-1}\ti u_{j'}(i))_{1\le j,j' \le n}\big)\big)^{-1}, 
\no \\ 
&\hspace*{8cm} \quad z\in \rho(A_1), \lb{2.29}
\end{align}
and (cf.~\eqref{2.23a})
\begin{equation}
\ti \theta_{2,j,j'}=(\ti u_j(i),\vartheta_2 \ti u_{j'}(i)), \quad 
1\leq j,j'\leq n. \lb{2.30}
\end{equation}
The proof of \eqref{2.28}--\eqref{2.30} is based on the following
elementary result.

\begin{lemma} \lb{2.3}
Let $\calH_n$, $n\in\bbN$ be an $n$-dimensional complex 
Hilbert space, $T\in\calB(\calH)$ a bounded linear operator 
in $\calH$ with $T^{-1}\in\calB(\calH)$. Assume that   
$\{\psi_j\}_{1\leq j\leq n}$ and  
$\{\ti \psi_j\}_{1\leq j\leq n}$ are {\rm(}not
necessarily  orthogonal\,{\rm)} bases in $\calH$. Then
\begin{align}
&\sum_{j,j'=1}^n \big(\big((\psi_\ell,T\psi_m)_{1\leq \ell,m\leq n}
\big)^{-1}\big)_{j,j'}(\psi_{j'},\cdot)\psi_j \no \\
&=\sum_{j,j'=1}^n \big(\big((\ti \psi_\ell,T\ti 
\psi_m)_{1\leq \ell,m\leq n}
\big)^{-1}\big)_{j,j'}(\ti \psi_{j'},\cdot)\ti \psi_j. \lb{2.31}
\end{align}
\end{lemma}

\section{The Direct Scattering Problem for Generalized \\Point 
Interactions} \lb{s3}

In the principal part of this section we apply the abstract 
framework surrounding Krein's formula \eqref{2.28} to the 
concrete situation of $n$ generalized point interactions  
in $\bbR^3$. At the end we derive the corresponding 
quantum mechanical scattering formalism, including explicit
expressions for the scattering wave functions and the 
scattering amplitude. 

In order to apply the results of Section~\ref{s2}, we now make a
series of identifications:
\begin{align}
& \calH=L^2(\bbR^3), \lb{3.1} \\
& A=\ol{-\Delta\big|_{C^\infty_0(\bbR^3\backslash
\{\xi_1,\dots,\xi_n\}})}\, , \quad 
\{\xi_1,\dots,\xi_n\}\subset\bbR^3, \,\, \xi_j\neq\xi_{j'} 
\text{ for } j\neq j', \lb{3.2} \\
&\ker(A^*-i)=\text{span}\{\ti u_j(i,x)=
G_0(i,x-\xi_j)\}_{1\leq j\leq n}, \lb{3.3} \\
&G_0(z,x-y)=(-\Delta -z)^{-1}(x,y)=
\f{\exp(iz^{1/2}|x-y|)}{4\pi |x-y|}, \lb{3.4} \\
& \hspace*{1.85cm} z\in\bbC\backslash\bbR, \,\, \Im(z^{1/2})>0, \,\,
x,y\in\bbR^3, \,  x\neq y, \, \no \\
& A_1=-\Delta, \quad \dom(-\Delta)=H^{2,2}(\bbR^3), \lb{3.5} \\
&\ti u_j(z)=(-\Delta-i)(-\Delta -z)^{-1}\ti u_j(i)=
G_0(z,\cdot -\xi_j), \lb{3.6} \\
& \hspace*{2.4cm} z\in\bbC\backslash\bbR, \, \Im(z^{1/2})>0,
\,\, 1\leq j\leq n. \no
\end{align}

In particular, a comparison of \eqref{3.3} and 
\begin{equation}
\ker(A^*-z)=\text{span}\{\ti u_j(z,x)=
G_0(z,x-\xi_j)\}_{1\leq j\leq n}, \quad z\in\bbC\backslash\bbR 
\lb{3.6a}
\end{equation}
shows that $\{\ti u_j(z)\}_{1\leq j\leq n}$ is a natural (though,  
not orthogonal) basis of $\ker(A^*-z)$. 
 
We note that the fact \eqref{3.3} can be found, for instance, 
in \cite[Sect.~II.1.1]{AGHKH88} and \cite{Zo80}.

Straightforward computations using 
\begin{equation}
(\pm i)^{1/2}=2^{-1/2}(\pm 1+i), \quad \ol{i(i)^{1/2}}=
i(-i)^{1/2} \lb{3.7}
\end{equation}
and the first resolvent equation
\begin{align}
&(-\Delta -z_1)^{-1}(-\Delta -z_2)^{-1}=(z_1 -z_2)^{-1}
[(-\Delta -z_1)^{-1}  - (-\Delta -z_2)^{-1}], \lb{3.8} \\
& \hspace*{8.5cm} z_1, z_2 \in \rho(-\Delta)  \no
\end{align}
repeatedly, then yield the following results.

\begin{lemma} \lb{l3.1} Let $z\in\rho(-\Delta)$ and 
$j,j'\in \{1,\dots,n\}$. Then
\begin{align}
& (\ti u_j(i),\ti u_{j'}(i))=\begin{cases}\|u_j(i)\|^2=(4\pi
2^{1/2})^{-1}, & j=j',\\
\Im(G_0(i,\xi_j-\xi_{j'})), & j \neq j', 
\end{cases} \lb{3.9} \\
& (\ti u_j(i),(-\Delta -z)^{-1}\ti u_j(i))=(4\pi)^{-1}
(z^2 +1)^{-1}[iz^{1/2}-i(-i)^{1/2}-2^{-1/2}(z+i)], 
\lb{3.10} \\
& (\ti u_j(i),(-\Delta -z)^{-1}\ti u_{j'}(i))=(z^2 +1)^{-1}
[G_0(z,\xi_j-\xi_{j'})-G_0(-i,\xi_j-\xi_{j'})] \no \\
&\hspace*{4.3cm} -(z-i)^{-1}\Im(G_0(i,\xi_j-\xi_{j'})), 
\quad j\neq j', \lb{3.11} \\
& \Im(z) ((\ti u_j(i),\ti u_{j'}(i))=\begin{cases} 
 (4\pi)^{-1}\Re (z^{1/2}), & j=j', \\
\Im(G_0(z,\xi_j-\xi_{j'})), & j\neq j', 
\end{cases} \lb{3.12}
\end{align}
where
\begin{equation}
\Im(G_0(z,\xi_j-\xi_{j'}))=[G_0(z,\xi-\xi_{j'})-
G_0(\ol z, \xi_j-\xi_{j'})]/(2i). \lb{3.12a}
\end{equation}
\end{lemma}

Given these preliminaries, one can now describe the $n^2$ 
(real) parameter family of all self-adjoint extensions of $A$, 
relatively prime to $A_1=-\Delta$ with respect to $A$, 
by appealing to Krein's formula \eqref{2.28}, \eqref{2.29} as 
follows. (In passing we note that $A_1=-\Delta$, as defined in 
\eqref{3.5}, is the Friedrichs extension of $A$.)~One defines
\begin{equation}
\Xi=\{\xi_1,\dots,\xi_n\} \subset \bbR^3 \lb{3.13}
\end{equation}
and denotes by 
\begin{equation}
\theta =\big(\theta_{j,j'}\big)_{1\leq j,j'\leq n} =\theta^* 
\lb{3.14}
\end{equation}
a self-adjoint $n\times n$ matrix in $\bbC^n$. Combining Krein's 
formula \eqref{2.28}, \eqref{2.29} with \eqref{3.10} and
\eqref{3.11} then yields the principal result of this section.

\begin{theorem} \lb{t3.2}
Let $z\in\bbC\backslash\bbR$. Then the $n^2$ {\rm(}real\,{\rm)}
parameter family of all self-adjoint  extensions
$-\Delta_{\theta,\Xi}$  of
$\ol{-\Delta\big|_{C^\infty_0(\bbR^3\backslash\Xi)}}$\,, 
relatively prime to $-\Delta$ with respect to 
$\ol{-\Delta\big|_{C^\infty_0(\bbR^3\backslash\Xi)}}$\,,    
can be parametrized by all self-adjoint 
$n\times n$ matrices $\theta$ in $\bbC^n$ with
\begin{equation}
\{((1/2) +m)\pi\}_{m\in\bbZ}\notin\spec(\theta). \lb{3.15}
\end{equation} 
An explicit representation for $-\Delta_{\theta,\Xi}$ is 
provided by  
\begin{equation}
(-\Delta_{\theta,\Xi} -z)^{-1}=(-\Delta -z)^{-1} + 
\sum_{j,j'=1}^n P_{\theta,\Xi}(z)_{j,j'}
(G_0(\ol z,\cdot -\xi_{j'}), \cdot )G_0(z,\cdot -\xi_j), 
\lb{3.16}
\end{equation}
where
\begin{align}
& \big(P_{\theta,\Xi}(z)^{-1}\big)_{j,j'} \lb{3.17} \\
&=\begin{cases} 
-(4\pi)^{-1}iz^{1/2}-(4\pi)^{-1}2^{-1/2}+ 
(\tan(\theta/2)_{j,j}, &j=j', \\
-G_0(z,\xi_j-\xi_{j'})+\Re(G_0(i,\xi_j-\xi_{j'})) + 
(\tan(\theta/2))_{j,j'}, &j\neq j', 
\end{cases} \no
\end{align}
and
\begin{equation}
\Re(G_0(z,\xi_j-\xi_{j'}))=[G_0(z,\xi_j-\xi_{j'})+
G_0(\ol z, \xi_j-\xi_{j'})]/2. \lb{3.18}
\end{equation}

\end{theorem}

\begin{remark} \lb{r3.3}
(i) Whenever $\theta$ in \eqref{3.16} has an eigenvalue 
$((1/2) +m_0)\pi$ for some $m_0\in\bbZ$, 
$\wti P_{\theta,\Xi}(z)$ becomes a singular matrix, 
$\det (\wti P_{\theta,\Xi}(z))=0$, $z\in\bbC\backslash\bbR$. 
In this case at least one point $\xi_{j_0}$ is removed {}from 
$\Xi$ and one effectively considers self-adjoint extensions 
of  
$A=\ol{-\Delta\big|_{C^\infty_0(\bbR^3\backslash
\{\Xi\backslash\{\xi_{j_0}\}\})}}$\,, 
parametrized in terms of $(n-1)\times (n-1)$ (or less) 
dimensional self-adjoint matrices $\theta$. In particular, 
the Friedrichs extension of 
$A=\ol{-\Delta\big|_{C^\infty_0(\bbR^3\backslash\Xi)}}$, 
given by $A_1=-\Delta$, formally corresponds to the extreme 
case $\theta =\pi I_n$ in \eqref{3.16}, \eqref{3.17}. \\
(ii) It seems appropriate to call the $n^2$-parameter family 
$-\Delta_{\theta,\Xi}$ defined by \eqref{3.16}, 
\eqref{3.17} the {\it
generalized} point interaction Hamiltonian, distinguishing it 
{}from the usually considered $n$-parameter family of (local) 
point interactions. In fact, introducing $\ul \alpha =
(\alpha_1,\dots,\alpha_n)\in\bbR^n$, the standard 
$n$-parameter family of self-adjoint extensions 
$-\Delta_{\ul\alpha,\Xi}$ of 
$\ol{-\Delta\big|_{C^\infty_0(\bbR^3\backslash\Xi)}}$ 
emerges as a special case of \eqref{3.16}, \eqref{3.17} by 
choosing
\begin{equation}
\tan(\theta/2)=\big((\alpha_j+(4\pi)^{-1}2^{-1/2})
\delta_{j,j'} - \Re(\wti G_0(i,\xi_j-\xi_{j'}))
\big)_{1\leq j,j'\leq n}, \lb{3.19}
\end{equation}
where
\begin{equation}
\wti G_0(z,x)=\begin{cases}G_0(z,x)&\text{if } x\neq 0, \\
0&\text{if } x=0 \text{ (or if $n=1$),} \end{cases} 
\lb{3.20}
\end{equation}
and $\Re(\wti G_0(z,\xi_j-\xi_{j'}))=[\wti G_0(z,\xi-\xi_{j'})
+\wti G_0(\ol z, \xi_j-\xi_{j'})]/2$. Insertion of 
\eqref{3.19} into \eqref{3.16}, \eqref{3.17} then yields
\begin{equation}
(-\Delta_{\ul\alpha,\Xi} -z)^{-1}=(-\Delta -z)^{-1} + 
\sum_{j,j'=1}^n \big(\Gamma_{\ul\alpha,\Xi}
(z)^{-1}\big)_{j,j'}
(G_0(\ol z,\cdot -\xi_{j'}), \cdot )G_0(z,\cdot -\xi_j), 
\lb{3.21}
\end{equation}
where
\begin{equation}
\Gamma_{\ul\alpha,\Xi}(z)
=\big((-(4\pi)^{-1}iz^{1/2}+\alpha_j)\delta_{j,j'}
-\wti G_0(z,\xi_j-\xi_{j'})\big)_{1\leq j,j'\leq n},  \lb{3.22}
\end{equation}
in accordance with \cite[p.~113]{AGHKH88}. While most efforts 
in connection with finitely many point interactions focus on 
the $n$-parameter family $-\Delta_{\ul\alpha,\Xi}$ (cf. the 
detailed discussion in \cite[Ch.~II.1]{AGHKH88} and the 
references therein), the general $n^2$-parameter family of 
generalized point interactions has been discussed by 
Dabrowski and Grosse \cite{DG85} in 1985. The treatment in 
\cite{DG85} also combines Krein's resolvent formula with von
Neumann's parametrization of self-adjoint extensions, but is 
somewhat less detailed than our present approach. (In 
particular, their matrix $S(\ol z,z_0)$, and hence their 
$M(z)$, are not explicitly computed in section~II of
\cite{DG85}, although these quantities can be inferred {}from the
scaling limit approach in section~IV via their formula 
(4.18).)    
\end{remark}

Finally, we briefly discuss stationary quantum scattering theory
following the lines of \cite[Sect.~II.1.5]{AGHKH88}
and \cite{Ra86}. Given the
resolvent kernel of $-\Delta_{\theta,\Xi}$ in \eqref{3.16}, one 
computes
\begin{align}
&\lim_{\varepsilon\downarrow 0}
\underset{\substack{ |y|\to\infty\\ -|y|^{-1}y=\omega}}
{\lim} 
4\pi |y|e^{-i(k+i\varepsilon)|y|}
(-\Delta_{\theta,\Xi} -(k+i\varepsilon)^2)^{-1}(x,y) \no \\
&=e^{ik\omega \cdot x}+\sum_{j,j'=1}^n P_{\theta,\Xi}(k^2)_{j,j'}
e^{ik\omega\cdot\xi_{j'}}G_0(k^2,x-\xi_j) \lb{3.24} \\
&=\Psi_{\theta,\Xi}(x,k,\omega), \quad k\in\bbR, \,\,
\det(P_{\theta,\Xi}(k^2))\neq 0, 
\,\, \omega\in S^2, \,\, x\in\bbR^3\backslash\Xi\,. \no
\end{align}
Moreover, since
\begin{equation}
(-\Delta\Psi)(x,k,\omega)=k^2\Psi(x,k,\omega), \quad 
x\in\bbR^3\backslash\Xi  \lb{3.25}
\end{equation}
in the distributional sense as well as pointwise,
$\Psi_{\theta,\Xi}(k\omega,x)$, 
$k\in\bbR$, $\det(P_{\theta,\Xi}(k^2))\neq 0$,  
$\omega\in S^2$, $x\in\bbR^3\backslash\Xi$, represent the 
generalized eigenfunctions, that is, the quantum scattering 
wave functions associated with $-\Delta_{\theta,\Xi}$. 

The corresponding quantum scattering amplitude 
$A_{\theta,\Xi}(\omega',\omega,k)$ is then computed as 
follows,
\begin{align}
A_{\theta,\Xi}(\omega',\omega,k)&=
\underset{\substack{ |x|\to\infty\\ |x|^{-1}x=\omega'}}
{\lim} |x|e^{-ik|x|}\big[\Psi_{\theta,\Xi}(x,k,\omega)
-e^{ik\omega\cdot x}\big] \no \\
&=(4\pi)^{-1}\sum_{j,j'=1}^n P_{\theta,\Xi}(k^2)_{j,j'}
e^{ik(\omega\cdot\xi_{j'}-\omega'\cdot\xi_j)}, \lb{3.26} \\
&\hspace*{.45cm} k\in\bbR, \,\, \det(P_{\theta,\Xi}(k^2))\neq 0, 
\,\, \omega, \omega'\in S^2. \no
\end{align}
The corresponding scattering matrix $S_{\theta,\Xi}(k)$ in 
$L^2(S^2)$ is then given by 
\begin{align}
&S_{\theta,\Xi}(k)=I+\f{ik}{8\pi^2}\sum_{j,j'=1}^n
P_{\theta,\Xi}(k^2)_{j,j'}\big(e^{-ik\xi_{j'}\cdot (\cdot)},
\cdot\big) e^{-ik\xi_j\cdot (\cdot)}, \lb{3.27} \\
&\hspace*{4.72cm} k\in\bbR, \,\, \det(P_{\theta,\Xi}(k^2))\neq 0. 
\no
\end{align}

\begin{remark} \lb{r3.4}
Since $S_{\theta,\Xi}(k)$ is unitary in $L^2(S^2)$ (this either 
follows {}from abstract methods since $-\Delta$ and
$-\Delta_{\theta,\Xi}$ are self-adjoint and the second term on the
right-hand side of
\eqref{3.16} is of rank $n$ and hence a  trace class operator, or
directly {}from \eqref{3.17} and 
\eqref{3.27}), the scattering amplitude
$A_{\theta,\Xi}(\omega',\omega,k)$ (the integral kernel of 
$S_{\theta,\Xi}(k)-I$) automatically satisfies the 
(generalized) {\it optical theorem}, that is, 
\begin{align}
&\Im(A_{\theta,\Xi}(\omega',\omega,k))=(4\pi)^{-1}k\int_{S^2} 
d\omega''\, A_{\theta,\Xi}(\omega'',\omega,k)
\ol{A_{\theta,\Xi}(\omega'',\omega',k)}\, , \lb{3.28} \\
&\hspace*{4.18cm} k\in\bbR, \,\, \det(P_{\theta,\Xi}(k^2))\neq 0, 
\,\, \omega,\omega'\in S^2. \no
\end{align}

On the other hand, {\it reciprocity} of the scattering 
amplitude $A_{\theta,\Xi}(\omega',\omega,k)$, defined by
\begin{equation}
A_{\theta,\Xi}(\omega',\omega,k)=
A_{\theta,\Xi}(-\omega,-\omega',k), \quad 
k\in\bbR, \,\, \det(P_{\theta,\Xi}(k^2))\neq 0, 
\,\, \omega,\omega' \in S^2 \lb{3.29}
\end{equation}
is satisfied if and only if
\begin{equation}
\theta^t=\theta, \lb{3.30}
\end{equation}
where $T^t$ denotes the transpose of the matrix $T$. Together 
with the requirement of self-adjointness of $\theta$, 
$\theta^*=\theta$, this yields an $n(n+1)/2$ (real) parameter 
family of operators $-\Delta_{\theta,\Xi}$ satisfying 
$\theta=\theta^*=\theta^t$. (The number of real elements 
above and on the diagonal of $\theta$ equals $\sum_{j=1}^n 
j =n(n+1)/2$.)

Similarly, the {\it reality} constraint on
$A_{\theta,\Xi}(\omega',\omega,k)$, that is, the requirement 
\begin{equation}
\ol{A_{\theta,\Xi}(\omega',\omega,k)}=
A_{\theta,\Xi}(\omega',\omega,-k), \quad 
k\in\bbR, \,\, \det(P_{\theta,\Xi}(k^2))\neq 0, 
\,\, \omega,\omega' \in S^2 \lb{3.31}
\end{equation}
is satisfied if and only if $\theta$ is a real matrix,
\begin{equation}
\theta_{j,j'}=\ol{\theta_{j,j'}}, \quad 1\leq j,j' \leq n. 
\lb{3.32}
\end{equation}
Together with self-adjointness of $\theta$ this again results 
in $\theta=\theta^*=\theta^t$ and hence is equivalent to the 
reciprocity requirement. (For background material on
properties of the scattering amplitude, such as the optical
theorem, reciprocity, and reality, we refer to 
\cite[Sect.~I.4]{Ra86}
for obstacle scattering and \cite[Sect.~3.6]{Th81} in the
context of potential scattering.)

It is interesting to observe that these natural requirements 
on the scattering amplitude, such as the optical theorem, 
reciprocity, and reality, are satisfied for an 
$n(n+1)/2$-parameter family of generalized point interactions 
(though, not for the full $n^2$-parameter family) and hence 
for a larger family than the usually considered 
$n$-parameter family of (local) point interactions 
$-\Delta_{\ul\alpha,\Xi}$.  
\end{remark}

We conclude this section with the following remark on  
space dimensions other than three (the interested reader can
find  many more details in \cite{AGHKH88}).

\begin{remark} \lb{r3.5}
All results of this section immediately extend to the case of 
two space dimensions replacing the Green's function \eqref{3.4} 
of the three-dimensional Laplacian by the corresponding
two-dimensional  Green's function
\begin{align}
&G_0(z,x-y)=(-\Delta -z)^{-1}(x,y)=(i/4)H_0^{(1)}
(z^{1/2}|x-y|), \lb{3.33} \\
& \hspace*{4.15cm} \Im(z^{1/2})>0, \,\, 
x,y\in\bbR^2, \, x\neq y. \no
\end{align}
Here $H_0^{(1)}(\cdot)$ denotes the Hankel function of order zero 
and first kind (cf.~\cite[Sect.~9.1]{AS72}). There are only 
minor changes required in \eqref{3.9}, \eqref{3.10}, and 
\eqref{3.11} due to the $\ln(z)$-behavior of \eqref{3.33} as 
$z\to 0$. Analogous results apply to the one-dimensional case 
using 
\begin{align}
&G_0(z,x-y)=(-\Delta -z)^{-1}(x,y)=(i/2)z^{-1/2}
e^{iz^{1/2}|x-y|}, \lb{3.34} \\
& \hspace*{5.17cm} \Im(z^{1/2})>0, \,\, x,y\in\bbR. \no
\end{align}
The one-dimensional case, however, differs {}from the 
two and three-dimensional cases since 
$B=\ol{-\f{d^2}{dx^2}\big|_{C_0^\infty(\bbR\backslash
\{\xi_1,\dots,\xi_n\})}}$ now has deficiency indices $(2n,2n)$ 
(as opposed to $(n,n)$ in two and three dimensions, 
cf.~\eqref{3.3}). Consequently, $B$ admits a $4n^2$ (real)
parameter family of self-adjoint extensions and hence 
additional types of (generalized) point interactions in 
dimension one. The proper definition of $A$ with 
\begin{equation}
\ker(A^*-z)=\text{span}\{G_0(z,\cdot -\xi_j)\}
_{1\leq j\leq n} \lb{3.35}
\end{equation}
in one dimension is given by
\begin{equation}
A=-\f{d^2}{dx^2}, \quad \dom (A)=\{g\in H^{2,2}(\bbR) \,|\, 
g(\xi_j)=0, \, 1\leq j\leq n \}. \lb{3.36}
\end{equation}

Finally, since
\begin{equation}
\ol{-\Delta\big|_{C_0^\infty(\bbR^d \backslash
\{\xi_1,\dots,\xi_n\})}}=-\Delta\big|_{H^{2,2}(\bbR^d)} 
\text{ for } d\geq 4 \lb{3.37}
\end{equation} 
(i.e., $-\Delta\big|_{C_0^\infty(\bbR^d \backslash
\{\xi_1,\dots,\xi_n\})}$ is essentially self-adjoint for 
$d\geq 4$), there are no (generalized) point interactions 
in four dimensions or higher.
\end{remark}

\section{A Uniqueness Result} \lb{s4}

Given the preparations in Section~\ref{s3}, the principal
purpose of our final Section~\ref{s4} is to provide a
solution of the inverse problem IP~1.2 formulated in
Section~\ref{s1}. More precisely, we will prove the
following uniqueness result (we freely use the notation
established in Sections~\ref{s1}--\ref{s3} throughout this 
section).

\begin{theorem} \lb{t4.1}
Let $k_0>0$ and assume that $\det(P_{\theta,\Xi}(k_0^2))\neq
0$. Then the data
$\{G_{\theta,\Xi}(k_0^2,x,y)\}_{\substack{x,y\in  P\\x\neq
y}}$ uniquely determine
$\Xi=\{\xi_1,\dots,\xi_n\}\subset\bbR^3_-$ and the
self-adjoint $n\times n$ matrix $\theta$ in $\bbC^n$.
\end{theorem}
\begin{proof}
Given the data $G_{\theta,\Xi}(k_0^2,x,y)$ for all $x,y\in
P$, $x\neq y$, $\det(P_{\theta,\Xi}(k_0^2))\neq 0$, the
Poisson-type formula,
\begin{align}
& w(s,y)=\int_Pd^2\sigma\,G_{\theta,\Xi}(k_0^2,(\sigma,0),y)
\f{\partial}{\partial
\sigma_3}\big[G_0(k_0^2,(\sigma,\sigma_3)-s) \no \\
& \hspace*{6cm} -
G_0(k_0^2,(\sigma,-\sigma_3)-s)\big]\big|_{\sigma_3 =0}, 
\lb{4.1} \\
& \hspace*{8.28cm} s\in\bbR^3_+, \,\, y\in P \no
\end{align}
(here $\sigma=(\sigma_1,\sigma_2)\in\bbR^2$) yields the 
solution of the problem
\begin{equation}
\begin{cases}
(-\bigtriangledown_x^2-k_0^2)w(x,y)=0, \quad x\in\bbR_+^3,\\
\lim_{|x|\to\infty}|x|\big[\f{\partial}{\partial|x|}w(x,y)
-ik_0w(x,y)\big]=0 \text{ uniformly in directions
$\omega=|x|^{-1}x$} \\
\hspace*{4.25cm} \text{and uniformly in $y$ for $y$ varying
in compact sets,} \\
w(x,y)\big|_{x_3=0}=G_{\theta,\Xi}(k_0^2,x,y)
\end{cases} \lb{4.2}
\end{equation}
for each fixed $y\in P$. In particular, $w(x,y)$ in
\eqref{4.1} represents
\begin{equation}
G_{\theta,\Xi}(k_0^2,x,y) \text{ for all } x\in\bbR_+^3\cup
P, \, y\in P, \, x\neq y. \lb{4.3}
\end{equation}
By symmetry of the Green's function
$G_{\theta,\Xi}(k_0^2,x,y)$ with respect to $x$ and $y$,
\begin{equation}
G_{\theta,\Xi}(k_0^2,x,y)=G_{\theta,\Xi}(k_0^2,y,x), \quad
x,y\in\bbR^3, \, x\neq y, \lb{4.4}
\end{equation}
we also determined
\begin{equation}
G_{\theta,\Xi}(k_0^2,x,y) \text{ for all }
x\in P, \, y\in\bbR_+^3\cup P, \, x\neq y. \lb{4.5}
\end{equation}
Moreover, using $G_{\theta,\Xi}(k_0^2,(\sigma,0),y)$ with 
$y\in\bbR_+^3\cup P$ (instead of $y\in P$) in \eqref{4.1}
then determines
\begin{equation}
G_{\theta,\Xi}(k_0^2,x,y) \text{ for all }
x,y \in\bbR_+^3\cup P, \, x\neq y. \lb{4.6}
\end{equation}
In other words, we managed to lift the data {}from $P$ to
$\bbR^3_+\cup P$.

Next, the explicit formula \eqref{3.16} for
$(-\Delta_{\theta,\Xi}-z)^{-1}$ yields
\begin{align}
&G_{\theta,\Xi}(z,x,y)=G_0(z,x,y)+\sum_{j,j'=1}^n 
P_{\theta,\Xi}(z)_{j,j'}G_0(z,x-\xi_j)G_)(z,y-\xi_{j'}), 
\lb{4.7} \\
& \hspace*{5cm} \det(P_{\theta,\Xi}(z))\neq 0, \,\, 
x,y\in\bbR^3\backslash\Xi, \, x\neq y, \no 
\end{align}
with $P_{\theta,\Xi}(z)$ defined in \eqref{3.17}. 
Hence one concludes
\begin{equation}
(-\bigtriangledown^2_x-z)G_{\theta,\Xi}(z,x,y)=0, \quad 
x,y\in\bbR^3\backslash\Xi, \, x\neq y. \lb{4.8}
\end{equation}
Thus, the data $\{G_{\theta,\Xi}(k_0^2,x,
y)\}_{\substack{x,y\in  P\\x\neq y}}$, 
$\det(P_{\theta,\Xi}(k_0^2))\neq 0$ uniquely determine 
\begin{equation}
G_{\theta,\Xi}(k_0^2,x,y) \text{ for all }
x,y \in\bbR^3\backslash\Xi, \, x\neq y \lb{4.9}
\end{equation}
by the unique continuation property \cite[Sect.~17.2]{Ho85} 
applied to \eqref{4.8}.
The singularity structure of \eqref{4.7} then determines
$\xi_1,\dots,\xi_n$ and hence $\Xi$. Similarly, taking 
$x\to\xi_j$ and $y\to\xi_{j'}$ independently, determines 
$P_{\theta,\Xi}(k_0^2)_{j,j'}$, $1\leq j,j'\leq n$, and
hence $\theta$. Thus, Theorem~\ref{t4.1} is proved.
\end{proof}


\end{document}